\documentclass[11pt,reqno]{amsart}
\usepackage{graphicx}
\usepackage{verbatim}
\usepackage{textcomp}
\usepackage{amssymb}
\usepackage{cite}
\usepackage{amsmath}
\usepackage{latexsym}
\usepackage{amscd}
\usepackage{amsthm}
\usepackage{mathrsfs}
\usepackage{bm}
\usepackage{url}
\usepackage{hyperref}
\usepackage{bookmark}
\allowdisplaybreaks[3]
\newtheorem{thm}{Theorem}[section]

\newtheorem{prop}[thm]{Proposition}

\theoremstyle{definition}

\theoremstyle{remark}
\newtheorem{rem}{Remark}[section]
\numberwithin{equation}{section}
\begin{document}
\title[Some inequalities between Laplacian eigenvalues]
{Some inequalities between Laplacian eigenvalues on Riemannian manifolds }

\author{Guangyue Huang }
\author{Xuerong Qi}
\address{Department of Mathematics, Henan Normal
University, Xinxiang 453007, P.R. China}

\email{hgy@henannu.edu.cn (G. Huang) }

\address{School of Mathematics and Statistics, Zhengzhou University, Zhengzhou 450001,
P.R. China}

\email{xrqi@zzu.edu.cn (X. Qi)}

\thanks{The research of the authors is
supported by NSFC(Nos. 11971153,  11671121).
}

\begin{abstract}
In this paper, we study a first Dirichlet eigenfunction of the weighted $p$-Laplacian on a bounded domain in a complete weighted Riemannian manifold. By constructing gradient estimates for a first eigenfunction, we obtain some relationships between weighted $p$-Laplacian first eigenvalues. As an immediate application, we also obtain some eigenvalue comparison results between the first Dirichlet eigenvalue of the weighted Laplacian, the first clamped plate eigenvalue and the first buckling eigenvalue.

\end{abstract}



\maketitle

\section{Introduction}

Let $(M^n, g)$ be an $n$-dimensional complete Riemannian manifold endowed with a weighted measure $d\mu=e^{-f}dv_{g}$,
where $f$ is a smooth real-valued function on $M$ and $dv_g$ is the Riemannian measure associated to the Riemannian metric $g$.
The triple $(M^n, g,d\mu)$ is also called a  complete metric measure space.
Denote by $\nabla$ and $\Delta$ the gradient operator and the Laplacian on $M$ with respect to the metric $g$, respectively.
It is well-known that
the $f$-Laplacian $\Delta_f$ is defined by
$$\Delta_f \varphi:=e^f{\rm div}(e^{-f}\nabla \varphi)=\Delta \varphi-\langle\nabla f, \nabla \varphi\rangle,\ \ \ \ \ \forall \ \varphi\in
C^\infty_0(M),$$
which is symmetric with respect to
the $L^2(M)$ inner product under the weighted measure
$d\mu$, that is, for any two functions $\varphi,\psi\in
C^\infty_0(M^n)$, it holds that
$$\int_{M}\varphi\Delta_f \psi\,d\mu=\int_{M}\psi\Delta_f \varphi\,d\mu=-\int_{M} \langle\nabla\varphi,\nabla \psi\rangle\,d\mu.$$
The $f$-Laplacian $\Delta_f$ is also called the drifting Laplacian, the weighted Laplacian or the
Witten-Laplacian. In particular, when $f$ is constant, the operator $\Delta_f$ becomes the Laplace operator $\Delta$.

The Bakry-\'{E}mery Ricci curvature (see \cite{Bakry85,li05,Wei09}) associated
to the $f$-Laplacian is given by
\begin{align}\label{0902-Int-1}
{\rm Ric}_f={\rm Ric}+\nabla^2f,
\end{align}
where ${\rm Ric}$ and $\nabla^2f$ denote the Ricci curvature of $M$ and the Hessian of $f$, respectively.
For a bounded domain $\Omega$ in $M^n$, the $f$-mean curvature of the boundary $\partial \Omega$ is defined by
\begin{align}\label{0902-Int-2}
H_f=H-\langle\nabla f,\nu\rangle,
\end{align}
where $\nu$ denotes the outward unit normal to $\partial \Omega$ and $H$ is the mean curvature (the trace of the second fundamental form) of $\partial \Omega$.

Recently, in order to find an ``universal" constant $c$ on a mean-convex bounded domain $\Omega$ of a complete Riemannian manifold (for the research in the Euclidean space, see \cite{Polya1952,Szego1954,Aviles1986,Payne1955,
Levine1986,Friedlander1991,CCWX2012} and the references therein) such that for any positive integer $k$,
$$\mu_{k+1}(\Omega)<c\lambda_k(\Omega),$$
where $\mu_{k+1}$ is the $(k+1)$-th Neumann eigenvalue and $\lambda_k$ denotes the $k$-th Dirichlet eigenvalue of the Laplacian,
Ilias and Shouman \cite{Ilias2020} studied the relationships between the first Dirichlet eigenvalue, the first clamped plate eigenvalue and the first buckling eigenvalue. The main method used in \cite{Ilias2020} is the gradient estimate with respect to a first Dirichlet eigenfunction of the weighted Laplacian $\Delta_{f}$.

Let $p\in(1,+\infty)$, for a function $\phi\in W^{1,p}$, the well-known weighted $p$-Laplacian is defined by
\begin{align}\label{0902-Int-3}
\Delta_{p,f}\phi:=e^{f}{\rm div}(e^{-f}|\nabla \phi|^{p-2}\nabla \phi),
\end{align}
which is understood in distribution sense.
In this paper, we consider the following Dirichlet eigenvalue problem (see \cite{Wangli2012,liHuang2020,Naber2014,Wangli2009}):
\begin{equation}
\left\{\begin{array}{l}\label{0909-Int-1}
\Delta_{p,f} u=-\lambda_{p,f}|u|^{p-2} u, \ \ {\rm in}\ \Omega,\\
u|_{\partial \Omega}=0.
\end{array}\right.
\end{equation}

It is known that \eqref{0909-Int-1} has a weak solution in $W^{1,p}_{0}(\Omega)$, the completion of the set $C^{\infty}_{0}(\Omega)$ of
smooth functions compactly supported on $\Omega$ under the Sobolev norm
$$ \| u \|_{1,p}=\bigg\{\int_{\Omega}
\big(|u|^{p}+|\nabla u|^{p}\big)d\mu\bigg\}^{\frac{1}{p}}.$$
Furthermore, the solution is positive in $\Omega$ and is unique up to a multiplicative constant. The solution is called a first eigenfunction of the weighted $p$-Laplace operator.
The first Dirichlet eigenvalue $\lambda_{1;p,f}$ of the weighted $p$-Laplacian can be characterized by
$$\lambda_{1;p,f}={\rm inf}\bigg\{\frac{\int_{\Omega}|\nabla \phi|^{p}d\mu}{\int_{\Omega}|\phi|^{p}d\mu} \ \bigg| \ \phi\in W^{1,p}_{0}(\Omega),\ \phi\not\equiv 0\bigg\}.$$

Inspired by \cite{Ilias2020}, we first construct gradient estimates with respect to a first Dirichlet eigenfunction of the weighted $p$-Laplacian. Using the gradient estimates for a first Dirichlet eigenfunction, we prove the following theorem:

\begin{thm}\label{1-thm1}
Let $(M^{n},g,d\mu)$ be an $n$-dimensional complete metric measure space and $\Omega$ be a bounded domain. Denote by $u$ the positive weak solution corresponding to the first eigenvalue $\lambda_{1;p,f}$. For $p>1$, if either
$p\neq2$ and
\begin{align}\label{1-thm-11}
{\rm Ric}_f-(p-1)(p-2)\frac{\nabla^2u}{u}\geq0,
\end{align}
or $p=2$ and ${\rm Ric}_f\geq0$
with $H_f\geq0$ on $\partial \Omega$,
then the first eigenvalue of the problem \eqref{0909-Int-1} satisfies
\begin{align}\label{3-thm-33}
\frac{(p^2+4p-2)-p\sqrt{p^2+8p-4}}{2(p+1)(2p-1)}
\lambda_{1;p,f}^2&\leq\lambda_{1;2p,f}\notag\\
&\leq\frac{(p^2+4p-2)+p\sqrt{p^2+8p-4}}{2(p+1)(2p-1)}\lambda_{1;p,f}^2.
\end{align}

\end{thm}

When $p=2$, $\Delta_{2,f}$ becomes $\Delta_{f}$.
Let $\Gamma_{1;f}$ be the first nonzero eigenvalue of the following clamped plate problem
\begin{equation}
\left\{\begin{array}{l}\label{0823-Sec-27}
\Delta_f^2 u=\Gamma_f\, u, \ \ {\rm in}\ \Omega,\\
u|_{\partial \Omega}=\frac{\partial u}{\partial\nu}|_{\partial \Omega}=0,
\end{array}\right.
\end{equation}
and $\Lambda_{1;f}$ be the first nonzero eigenvalue of the following buckling problem
\begin{equation}
\left\{\begin{array}{l}\label{0823-Sec-23}
\Delta_f^2 u=-\Lambda_f  \Delta_f u, \ \ {\rm in}\ \Omega,\\
u|_{\partial \Omega}=\frac{\partial u}{\partial\nu}|_{\partial \Omega}=0.
\end{array}\right.
\end{equation}
Denote by $\lambda_{1;f}$ the first Dirichlet eigenvalue of the $f$-Laplacian $\Delta_{f}$, that is, $\lambda_{1;f}=\lambda_{1;2,f}$. Then, we have the following

\begin{thm}\label{2-thm2}
Let $(M^{n},g,d\mu)$ be an $n$-dimensional complete metric measure space and $\Omega$ be a bounded domain.

(1) If
\begin{align}\label{2-thm-22}
{\rm Ric}_f-\frac{4}{3}\frac{\nabla^2u}{u}\geq0,
\end{align}
where $u$ is a first positive Dirichlet eigenfunction of the $f$-Laplacian, then we have
\begin{align}\label{add}
\Gamma_{1;f}\leq&\frac{64}{45}\lambda_{1;f}^2,
\qquad
\Lambda_{1;f}\leq \frac{64}{45}\lambda_{1;f}.
\end{align}

(2) If
${\rm Ric}_f\geq0$
and $H_f\geq0$ on $\partial \Omega$, then we have
\begin{align}\label{add2-thm-222222}
\Gamma_{1;f}\leq&\frac{16}{3}\lambda_{1;f}^2,
\qquad
\Lambda_{1;f}\leq \frac{16}{3}\lambda_{1;f}.
\end{align}

\end{thm}

\begin{rem}
The authors in \cite{Ilias2020} obtain the case of $p=2$ for Theorem \ref{1-thm1}. The case of $p\neq2$ in Theorem \ref{1-thm1} is new. On the other hand, the authors in \cite{Ilias2020} obtain $\Gamma_{1;f}<\frac{16}{3}\lambda_{1;f}^2$ and $\Lambda_{1;f}< 4\lambda_{1;f}$ (see Theorem 3.3 in \cite{Ilias2020}). Obviously, our results \eqref{add} in Theorem \ref{2-thm2} are new and sharper than those.
\end{rem}

\section{Gradient estimates for a first Dirichlet eigenfunction}

In this section, we construct the gradient estimate with respect to a first Dirichlet eigenfunction associated to the weighted $p$-Laplacian, which plays an important role in the proof of the main results.

Let $u$ be a weak solution associated with $\lambda_{1;p,f}$ to the problem \eqref{0909-Int-1}, that is,
\begin{equation}\label{0909-Prof-1}
\left \{ \aligned
-\Delta_{p,f} u&=\lambda_{1;p,f}|u|^{p-2} u,\quad \text{${\rm in} \ \ \Omega;$} \\
u&>0,\quad \text{${\rm in} \ \ \Omega;$} \\
u&=0, \quad \text{${\rm on} \ \ \partial\Omega.$}
\endaligned \right.
\end{equation}
Multiplying both sides of the above equality with $|u|^{\alpha+1}$ and using integration by parts, where $\alpha$ is a positive constant to be determined, we have
\begin{align}\label{0909-Prof-2}
\lambda_{1;p,f}\int
_{\Omega}|u|^{\alpha+p}\,d\mu
&=-\int_{\Omega}|u|^{\alpha+1}\Delta_{p,f} u\,d\mu\notag\\
&=\int_{\Omega}|\nabla u|^{p-2}\langle\nabla u,\nabla |u|^{\alpha+1}\rangle\,d\mu\notag\\
&=(\alpha+1)\int_{\Omega} |u|^{\alpha}|\nabla u|^{p} \,d\mu,
\end{align}
which is equivalent to
\begin{align}\label{0909-Prof-3}
\int_{\Omega} |u|^{\alpha}|\nabla u|^{p} \,d\mu=\frac{1}{\alpha+1}\lambda_{1;p,f}\int_{\Omega}|u|^{\alpha+p}\,d\mu.
\end{align}
By a direct calculation, we obtain
\begin{align}\label{0909-Prof-4}
\lambda_{1;p,f}\int_{\Omega} |u|^{\alpha}|\nabla u|^{p} \,d\mu
&=-\int_{\Omega} |u|^{\alpha-p+1}|\nabla u|^{p}\Delta_{p,f}u  \,d\mu\notag\\
&=\int_{\Omega} |\nabla u|^{p-2}\langle\nabla u,\nabla(|u|^{\alpha-p+1}|\nabla u|^{p})\rangle\,d\mu\notag\\
&=(\alpha-p+1)\int_{\Omega} |u|^{\alpha-p}|\nabla u|^{2p}\,d\mu\notag\\
& \ \ \ +p\int_{\Omega} |u|^{\alpha-p+1}|\nabla u|^{2p-3}\langle\nabla u,\nabla|\nabla u|\rangle\,d\mu.
\end{align}

Let
\begin{align}\label{0909-Prof-5}
I_{p,\alpha}=\frac{\int_{\Omega} |u|^{\alpha-p}|\nabla u|^{2p}\,d\mu}{\int_{\Omega} |u|^{\alpha}|\nabla u|^{p} \,d\mu}.
\end{align}
Then, \eqref{0909-Prof-4} yields
\begin{align}\label{0909-Prof-6}
\lambda_{1;p,f}\int_{\Omega} |u|^{\alpha}|\nabla u|^{p} \,d\mu=&(\alpha-p+1)I_{p,\alpha}\int_{\Omega} |u|^{\alpha}|\nabla u|^{p} \,d\mu\notag\\
&+p\int_{\Omega} |u|^{\alpha-p+1}|\nabla u|^{2p-3}\langle\nabla u,\nabla|\nabla u|\rangle\,d\mu,
\end{align}
which is equivalent to
\begin{align}\label{0909-Prof-7}
&(\alpha-p+1)I_{p,\alpha}\int_{\Omega} |u|^{\alpha}|\nabla u|^{p}\,d\mu\notag\\
&=\int_{\Omega} \Big(\lambda_{1;p,f}|u|^{\frac{\alpha}{2}}|\nabla u|^{\frac{p}{2}}-p|u|^{\frac{\alpha}{2}-p+1}|\nabla u|^{\frac{3}{2}p-3}\langle\nabla u,\nabla|\nabla u|\rangle\Big)|u|^{\frac{\alpha}{2}}|\nabla u|^{\frac{p}{2}}\,d\mu.
\end{align}
It follows from the H\"{o}lder inequality that
\begin{align}\label{0909-Prof-8}
&(\alpha-p+1)^2 I_{p,\alpha}^2\Big(\int_{\Omega} |u|^{\alpha}|\nabla u|^{p}\,d\mu\Big)^2\notag\\
&\leq\int_{\Omega} \Big(\lambda_{1;p,f}|u|^{\frac{\alpha}{2}}|\nabla u|^{\frac{p}{2}}-p|u|^{\frac{\alpha}{2}-p+1}|\nabla u|^{\frac{3}{2}p-3}\langle\nabla u,\nabla|\nabla u|\rangle\Big)^2\,d\mu
\int_{\Omega}|u|^{\alpha}|\nabla u|^{p}\,d\mu.
\end{align}
Hence
\begin{align}\label{0909-Prof-9}
&(\alpha-p+1)^2 I_{p,\alpha}^2
\int_{\Omega} |u|^{\alpha}|\nabla u|^{p}\,d\mu\notag\\
\leq&\int_{\Omega} \Big(\lambda_{1;p,f}|u|^{\frac{\alpha}{2}}|\nabla u|^{\frac{p}{2}}-p|u|^{\frac{\alpha}{2}-p+1}|\nabla u|^{\frac{3}{2}p-3}\langle\nabla u,\nabla|\nabla u|\rangle\Big)^2\,d\mu\notag\\
=&\lambda_{1;p,f}^2\int_{\Omega}|u|^{\alpha}|\nabla u|^{p}\,d\mu+p^2\int_{\Omega}|u|^{\alpha-2p+2}|\nabla u|^{3p-6}\langle\nabla u,\nabla|\nabla u|\rangle^2\,d\mu\notag\\
&-2p\lambda_{1;p,f}\int_{\Omega}|u|^{\alpha-p+1}|\nabla u|^{2p-3}\langle\nabla u,\nabla|\nabla u|\rangle\,d\mu.
\end{align}

Using the divergence theorem again, we obtain
\begin{align}\label{addadd0909-Prof-11}
&-2p\lambda_{1;p,f}\int_{\Omega}|u|^{\alpha-p+1}|\nabla u|^{2p-3}\langle\nabla u,\nabla|\nabla u|\rangle\,d\mu\notag\\
=&-\frac{p}{p-1}\lambda_{1;p,f}\int_{\Omega}|u|^{\alpha-p+1}
\langle\nabla u,\nabla|\nabla u|^{2p-2}\rangle\,d\mu\notag\\
=&\frac{p}{p-1}\lambda_{1;p,f}\int_{\Omega}
\big[(\alpha-p+1)|u|^{\alpha-p}|\nabla u|^{2}+|u|^{\alpha-p+1}\Delta_{f}u\big]|\nabla u|^{2p-2}\,d\mu\notag\\
=&\frac{p(\alpha-p+1)}{p-1}\lambda_{1;p,f}
\int_{\Omega}|u|^{\alpha-p}|\nabla u|^{2p}\,d\mu\notag\\
&+\frac{p}{p-1}\lambda_{1;p,f}
\int_{\Omega}|u|^{\alpha-p+1}|\nabla u|^{p}\Big(-\lambda_{1;p,f}|u|^{p-2}u-\langle\nabla u,\nabla|\nabla u|^{p-2}\rangle\Big)\,d\mu\notag\\
=&\frac{p(\alpha-p+1)}{p-1}\lambda_{1;p,f}
\int_{\Omega}|u|^{\alpha-p}|\nabla u|^{2p}\,d\mu-\frac{p}{p-1}
\lambda_{1;p,f}^2\int_{\Omega}|u|^{\alpha}|\nabla u|^{p}\,d\mu\notag\\
&-\frac{p(p-2)}{p-1}\lambda_{1;p,f}
\int_{\Omega}|u|^{\alpha-p+1}|\nabla u|^{2p-3}\langle\nabla u,\nabla|\nabla u|\rangle\,d\mu,
\end{align}
which shows
\begin{align}\label{0909-Prof-11}
&-2p\lambda_{1;p,f}\int_{\Omega}|u|^{\alpha-p+1}|\nabla u|^{2p-3}\langle\nabla u,\nabla|\nabla u|\rangle\,d\mu\notag\\
=&2(\alpha-p+1)\lambda_{1;p,f}\int_{\Omega}|u|^{\alpha-p}|\nabla u|^{2p}\,d\mu-2\lambda_{1;p,f}^2\int_{\Omega}|u|^{\alpha}|\nabla u|^{p}\,d\mu.
\end{align}
Inserting \eqref{0909-Prof-11} into \eqref{0909-Prof-9} gives
\begin{align}\label{0909-Prof-12}
&\left[(\alpha-p+1)^2 I_{p,\alpha}^2+\lambda_{1;p,f}^2\right]\int_{\Omega} |u|^{\alpha}|\nabla u|^{p}\,d\mu\notag\\
\leq& p^2\int_{\Omega}|u|^{\alpha-2p+2}|\nabla u|^{3p-6}\langle\nabla u,\nabla|\nabla u|\rangle^2\,d\mu\notag\\
&+2(\alpha-p+1)\lambda_{1;p,f}\int_{\Omega}
|u|^{\alpha-p}|\nabla u|^{2p}\,d\mu.
\end{align}

In order to prove Theorem \ref{1-thm1}, we need the following gradient estimate with respect to a first Dirichlet eigenfunction associated to the weighted $p$-Laplacian:

\begin{prop}\label{2-prop1}
Let $(M^{n},g, d\mu)$ be a complete metric measure space and $\Omega$ be a bounded domain. Denote by $u$ the first eigenfunction corresponding to the first Dirichlet eigenvalue $\lambda_{1;p,f}$ of the weighted $p$-Laplacian with $p>1$.
If either $\alpha\neq2(p-1)$ and
\begin{align}\label{1-prop-11}
{\rm Ric}_f+(p-1)(\alpha-2p+2)\frac{\nabla^2u}{u}\geq0,
\end{align}
or $\alpha=2(p-1)$ and ${\rm Ric}_f\geq0$ with $H_f\geq0$ on $\partial \Omega$,
then we have
\begin{align}\label{1-prop-111}
M_1\leq\frac{I_{p,\alpha}}{\lambda_{1;p,f}}\leq M_2,
\end{align}
where
$$M_1=\frac{\Big[2(\alpha-p+1)+\frac{p^2}{2p-1}\Big]
-\frac{p}{\sqrt{2p-1}}\sqrt{4(\alpha-p+1)+\frac{p^2}{2p-1}}}
{2(\alpha-p+1)^2},$$
$$M_2=\frac{\Big[2(\alpha-p+1)
+\frac{p^2}{2p-1}\Big]+\frac{p}{\sqrt{2p-1}}
\sqrt{4(\alpha-p+1)+\frac{p^2}{2p-1}}}{2(\alpha-p+1)^2}.$$
\end{prop}

\proof
Recall that the weighted Bochner formula (see the formula (2.3) in \cite{Wangli2012}) with respect to $\Delta_{p,f}$ is given by
\begin{align}\label{0909-Prof-13}
\frac{p-1}{p}\Delta_{p,f}|\nabla u|^{p}=&|\nabla u|^{2p-4}\big[|{\rm Hess}\,u|_A^2+{\rm Ric}_f(\nabla u,\nabla u)\big]\notag\\
&+|\nabla u|^{p-2}\langle\nabla u,\nabla\Delta_{p,f}u\rangle,
\end{align}
where
$$|{\rm Hess}\,u|_A^2=|\nabla^{2}u|^2+2(p-2)|\nabla |\nabla u||^{2}
+(p-2)^{2}\frac{\langle \nabla u,\nabla |\nabla u|\rangle^{2}}{|\nabla u|^{2}}.$$
Thus, we have
\begin{align}\label{0909-Prof-17}
&\int_{\Omega}|u|^{\alpha-2p+2}|\nabla u|^{3p-4}|{\rm Hess}\,u|_A^2\,d\mu\notag\\
=&\frac{p-1}{p}\int_{\Omega} |u|^{\alpha-2p+2}|\nabla u|^{p}\Delta_{p,f}(|\nabla u|^{p})\,d\mu\notag\\
& -\int_{\Omega}|u|^{\alpha-2p+2}|\nabla u|^{2p-2}\langle\nabla u,\nabla\Delta_{p,f}u\rangle\,d\mu\notag\\
&-\int_{\Omega}|u|^{\alpha-2p+2}|\nabla u|^{3p-4}{\rm Ric}_f(\nabla u,\nabla u)\,d\mu.
\end{align}

Next, we consider the following two cases:

{\it Case I.}
If $\alpha\neq2(p-1)$, then we also have
\begin{align}\label{0909-Prof-18}
&\int_{\Omega} |u|^{\alpha-2p+2}|\nabla u|^{p}\Delta_{p,f}(|\nabla u|^{p})\,d\mu\notag\\
=&-p(\alpha-2p+2)\int_{\Omega}|u|^{\alpha-2p+1}|\nabla u|^{3p-3}\langle\nabla u,\nabla|\nabla u|\rangle\,d\mu\notag\\
&-p^2\int_{\Omega}|u|^{\alpha-2p+2}|\nabla u|^{3p-4}|\nabla|\nabla u||^2\,d\mu
\end{align}
and
\begin{align}\label{0909-Prof-19}
&\int_{\Omega}|u|^{\alpha-2p+2}|\nabla u|^{2p-2}\langle\nabla u,\nabla\Delta_{p,f}u\rangle\,d\mu\notag\\
=&-(p-1)\lambda_{1;p,f}\int_{\Omega}|u|^{\alpha-p}|\nabla u|^{2p}\,d\mu.
\end{align}
Putting \eqref{0909-Prof-18} and \eqref{0909-Prof-19} into \eqref{0909-Prof-17} yields
\begin{align}\label{0909-Prof-20}
&\int_{\Omega}|u|^{\alpha-2p+2}|\nabla u|^{3p-4}|{\rm Hess}\,u|_A^2\,d\mu\notag\\
=&-\int_{\Omega}|u|^{\alpha-2p+2}|\nabla u|^{3p-4}{\rm Ric}_f(\nabla u,\nabla u)\,d\mu\notag\\
&+(p-1)\lambda_{1;p,f}\int_{\Omega}|u|^{\alpha-p}|\nabla u|^{2p}\,d\mu\notag\\
&-(p-1)(\alpha-2p+2)\int_{\Omega}|u|^{\alpha-2p+1}|\nabla u|^{3p-3}\langle\nabla u,\nabla|\nabla u|\rangle\,d\mu\notag\\
&-p(p-1)\int_{\Omega}|u|^{\alpha-2p+2}|\nabla u|^{3p-4}|\nabla|\nabla u||^2\,d\mu.
\end{align}

For any point $p\in \Omega$, we choose an orthonormal frame $\{e_i\}_{i=1}^n$ such that $\nabla u\parallel e_n$. Then, we have $\nabla u=|\nabla u|e_n$ and hence $u_n=|\nabla u|$, $u_1=\cdots=u_{n-1}=0$. At the point $p$, we have
\begin{align}\label{0909-Prof-14}
\langle\nabla u,\nabla|\nabla u|\rangle=u_nu_{nn}
\end{align}
and
\begin{align}\label{0909-Prof-15}
|{\rm Hess}\,u|_A^2&=\sum_{i,j}u_{ij}^2+2(p-2)\sum_{k=1}^{n}u_{kn}^2
+(p-2)^{2}u_{nn}^2\notag\\
&=(p-1)^2u_{nn}^2+2(p-1)\sum_{k=1}^{n-1}u_{kn}^2
+\sum_{k,l=1}^{n-1}u_{kl}^2\notag\\
&\geq(p-1)^2u_{nn}^2.
\end{align}
It follows from \eqref{0909-Prof-14} and \eqref{0909-Prof-15} that
\begin{align}\label{0909-Prof-16}
|\nabla u|^2|{\rm Hess}\,u|_A^2\geq (p-1)^2\langle\nabla u,\nabla|\nabla u|\rangle^2.
\end{align}
Using the Cauchy inequality, it holds that
\begin{align}\label{0909-Prof-21}
\langle\nabla u,\nabla|\nabla u|\rangle^2\leq|\nabla u|^{2}|\nabla|\nabla u||^2.
\end{align}
Inserting \eqref{0909-Prof-16} and \eqref{0909-Prof-21} into \eqref{0909-Prof-20} gives
\begin{align}\label{0909-Prof-22}
&(p-1)^2\int_{\Omega}|u|^{\alpha-2p+2}|\nabla u|^{3p-6}\langle\nabla u,\nabla|\nabla u|\rangle^2\,d\mu\notag\\
\leq&\int_{\Omega}|u|^{\alpha-2p+2}|\nabla u|^{3p-4}|{\rm Hess}\,u|_A^2\,d\mu\notag\\
\leq&-\int_{\Omega}|u|^{\alpha-2p+2}|\nabla u|^{3p-4}{\rm Ric}_f(\nabla u,\nabla u)\,d\mu\notag\\
&+(p-1)\lambda_{1;p,f}\int_{\Omega}|u|^{\alpha-p}|\nabla u|^{2p}\,d\mu\notag\\
&-(p-1)(\alpha-2p+2)\int_{\Omega}|u|^{\alpha-2p+1}|\nabla u|^{3p-3}\langle\nabla u,\nabla|\nabla u|\rangle\,d\mu\notag\\
&-p(p-1)\int_{\Omega}|u|^{\alpha-2p+2}|\nabla u|^{3p-6}\langle\nabla u,\nabla|\nabla u|\rangle^2\,d\mu,
\end{align}
which shows that
\begin{align}\label{0909-Prof-23}
&(p-1)(2p-1)\int_{\Omega}|u|^{\alpha-2p+2}|\nabla u|^{3p-6}\langle\nabla u,\nabla|\nabla u|\rangle^2\,d\mu\notag\\
\leq&-\int_{\Omega}|u|^{\alpha-2p+2}|\nabla u|^{3p-4}{\rm Ric}_f(\nabla u,\nabla u)\,d\mu+(p-1)\lambda_{1;p,f}\int_{\Omega}|u|^{\alpha-p}|\nabla u|^{2p}\,d\mu\notag\\
&-(p-1)(\alpha-2p+2)\int_{\Omega}|u|^{\alpha-2p+1}|\nabla u|^{3p-3}\langle\nabla u,\nabla|\nabla u|\rangle\,d\mu,
\end{align}
that is,
\begin{align}\label{0909-Prof-24}
&p^2\int_{\Omega}|u|^{\alpha-2p+2}|\nabla u|^{3p-6}\langle\nabla u,\nabla|\nabla u|\rangle^2\,d\mu\notag\\
\leq&-\frac{p^2}{(p-1)(2p-1)}\int_{\Omega}|u|^{\alpha-2p+2}|\nabla u|^{3p-4}{\rm Ric}_f(\nabla u,\nabla u)\,d\mu\notag\\
&-\frac{p^2(\alpha-2p+2)}{2p-1}\int_{\Omega}|u|^{\alpha-2p+1}|\nabla u|^{3p-3}\langle\nabla u,\nabla|\nabla u|\rangle\,d\mu\notag\\
&+\frac{p^2}{2p-1}\lambda_{1;p,f}\int_{\Omega}|u|^{\alpha-p}|\nabla u|^{2p}\,d\mu.
\end{align}
Therefore, putting \eqref{0909-Prof-24} into \eqref{0909-Prof-12} gives the estimate
\begin{align}\label{0909-Prof-25}
&(\alpha-p+1)^2 I_{p,\alpha}^2\int_{\Omega} |u|^{\alpha}|\nabla u|^{p}\,d\mu\notag\\
\leq&-\lambda_{1;p,f}^2\int_{\Omega}|u|^{\alpha}|\nabla u|^{p}\,d\mu\notag\\
&-\frac{p^2}{(p-1)(2p-1)}\int_{\Omega}|u|^{\alpha-2p+2}|\nabla u|^{3p-4}{\rm Ric}_f(\nabla u,\nabla u)\,d\mu\notag\\
&-\frac{p^2(\alpha-2p+2)}{2p-1}\int_{\Omega}|u|^{\alpha-2p+2}|\nabla u|^{3p-4}\frac{\nabla^{2}u(\nabla u,\nabla u)}{u}\,d\mu\notag\\
&+\Big[2(\alpha-p+1)+\frac{p^2}{2p-1}\Big]
\lambda_{1;p,f}\int_{\Omega}|u|^{\alpha-p}|\nabla u|^{2p}\,d\mu\notag\\
\leq&-\lambda_{1;p,f}^2\int_{\Omega}|u|^{\alpha}|\nabla u|^{p}\,d\mu\notag\\
&+\Big[2(\alpha-p+1)+\frac{p^2}{2p-1}\Big]
\lambda_{1;p,f}\int_{\Omega}|u|^{\alpha-p}|\nabla u|^{2p}\,d\mu,
\end{align}
where the last inequality we use \eqref{1-prop-11}. In particular, \eqref{0909-Prof-25} is equivalent to
\begin{align}\label{0909-Prof-26}
\lambda_{1;p,f}^2-\Big[2(\alpha-p+1)+\frac{p^2}{2p-1}\Big]\lambda_{1;p,f}I_{p,\alpha}+
(\alpha-p+1)^2 I_{p,\alpha}^2\leq0.
\end{align}
Solving this quadratic inequality with respect to $\frac{I_{p,\alpha}}{\lambda_{1;p,f}}$, it yields the desired estimate \eqref{1-prop-111}.

{\it Case II.}
If $\alpha=2(p-1)$, then \eqref{0909-Prof-17} becomes
\begin{align}\label{0914-Prof-1}
&\int_{\Omega}|\nabla u|^{3p-4}|{\rm Hess}\,u|_A^2\,d\mu\notag\\
=&\frac{p-1}{p}\int_{\Omega} |\nabla u|^{p}\Delta_{p,f}(|\nabla u|^{p})\,d\mu-\int_{\Omega}|\nabla u|^{3p-4}{\rm Ric}_f(\nabla u,\nabla u)\,d\mu\notag\\
&-\int_{\Omega}|\nabla u|^{2p-2}\langle\nabla u,\nabla\Delta_{p,f}u\rangle\,d\mu.
\end{align}

Furthermore,
\begin{align}\label{0914-Prof-3}
\int_{\Omega}|\nabla u|^{2p-2}\langle\nabla u,\nabla\Delta_{p,f}u\rangle\,d\mu
=-(p-1)\lambda_{1;p,f}\int_{\Omega}|u|^{p-2}|\nabla u|^{2p}\,d\mu,
\end{align}
\begin{align}\label{0914-Prof-2}
&\int_{\Omega} |\nabla u|^{p}\Delta_{p,f}(|\nabla u|^{p})\,d\mu\notag\\
=&-p^{2}\int_{\Omega} |\nabla u|^{3p-4}|\nabla|\nabla u||^2\,d\mu
+\int_{\partial \Omega} |\nabla u|^{2p-2}\langle\nu,\nabla|\nabla u|^{p}\rangle\,d\sigma\notag\\
=&-p^{2}\int_{\Omega} |\nabla u|^{3p-4}|\nabla|\nabla u||^2\,d\mu
-p\int_{\partial \Omega} |\nabla u|^{3p-2}H_f\,d\sigma\notag\\
&\leq-p^{2}\int_{\Omega} |\nabla u|^{3p-4}|\nabla|\nabla u||^2\,d\mu,
\end{align}
where we used $H_f\geq0$ and
$$\frac{1}{2}\frac{\partial(|\nabla u|^2)}{\partial \nu}=\frac{\partial u}{\partial \nu}\nabla^{2}u(\nu,\nu)=-H_f|\nabla u|^2$$
on the boundary $\partial \Omega$. Here we denote by $d\sigma$ the weighted measure induced on $\partial \Omega$.

Applying \eqref{0914-Prof-3} and \eqref{0914-Prof-2} into \eqref{0914-Prof-1}, we obtain
\begin{align}\label{0914-Prof-4}
&\int_{\Omega}|\nabla u|^{3p-4}|{\rm Hess}\,u|_A^2\,d\mu\notag\\
\leq&-\int_{\Omega}|\nabla u|^{3p-4}{\rm Ric}_f(\nabla u,\nabla u)\,d\mu+(p-1)\lambda_{1;p,f}\int_{\Omega}|u|^{p-2}|\nabla u|^{2p}\,d\mu\notag\\
&-p(p-1)\int_{\Omega}|\nabla u|^{3p-4}|\nabla|\nabla u||^2\,d\mu.
\end{align}
Therefore, following the proof as in the case of $\alpha\neq2(p-1)$, we derive the estimate \eqref{0909-Prof-26}
and the proof of Proposition \ref{2-prop1} is finished.

\section{Proof of theorems}

In this section, we prove the main results, namely Theorem 1.1 and Theorem 1.2.

\subsection{Proof of Theorem 1.1}
It is easy to see that \eqref{1-prop-111} is equivalent to
\begin{align}\label{0915-Prof-1}
M_1\lambda_{1;p,f}^2\leq\lambda_{1;p,f}I_{p,\alpha}\leq M_2\lambda_{1;p,f}^2.
\end{align}
Using \eqref{0909-Prof-3}, we have
\begin{align}\label{0915-Prof-2}
\lambda_{1;p,f}I_{p,\alpha}=(\alpha+1)\frac{\int_{\Omega} |u|^{\alpha-p}|\nabla u|^{2p}\,d\mu}{\int_{\Omega} |u|^{\alpha+p}\,d\mu},
\end{align}
which shows that \eqref{0915-Prof-1} is equivalent to
\begin{align}\label{0915-Prof-3}
\frac{1}{\alpha+1}M_1\lambda_{1;p,f}^2\leq
\frac{\int_{\Omega} |u|^{\alpha-p}|\nabla u|^{2p}\,d\mu}{\int_{\Omega} |u|^{\alpha+p}\,d\mu}\leq \frac{1}{\alpha+1} M_2\lambda_{1;p,f}^2.
\end{align}
In particular, choosing $\alpha=p$ in \eqref{0915-Prof-3} gives
\begin{align}\label{0915-Prof-4}
\frac{(p^2+4p-2)-p\sqrt{p^2+8p-4}}{2(p+1)(2p-1)}
\lambda_{1;p,f}^2&\leq\frac{\int_{\Omega}|\nabla u|^{2p}\,d\mu}{\int_{\Omega} |u|^{2p}\,d\mu}\notag\\
&\leq\frac{(p^2+4p-2)+p\sqrt{p^2+8p-4}}{2(p+1)(2p-1)}\lambda_{1;p,f}^2,
\end{align}
which shows that the estimate in Theorem 1.1 holds.

\subsection{ Proof of Theorem \ref{2-thm2}}
When $p=\alpha=2$, from Proposition \ref{2-prop1}, we know that if $H_f\geq0$ and ${\rm Ric}_f\geq0$, then
\begin{align}\label{0916-Prof-1}
\frac{1}{3}\leq\frac{I_{2,2}}{\lambda_{1;f}}
\leq3,
\end{align}
where
\begin{align}\label{0916-Prof-2}
I_{2,2}=\frac{\int_{\Omega} |\nabla u|^{4}\,d\mu}{\int_{\Omega} |u|^{2}|\nabla u|^{2} \,d\mu}.
\end{align}

Let $\varphi=u^{2}$, then the function $\varphi$ satisfies
$$\varphi|_{\partial \Omega}=\frac{\partial \varphi}{\partial \nu}\Big|_{\partial \Omega}=0.$$
Therefore,
\begin{align}\label{0916-Prof-3}
\Gamma_{1;f}&\leq \frac{\int_{\Omega}(\Delta_f\varphi)^2\,d\mu}
{\int_{\Omega}\varphi^2\,d\mu}\notag\\
&=\frac{4\int_{\Omega}\Big[u^2(\Delta_fu)^2+2u|\nabla u|^2\Delta_fu+|\nabla u|^4\Big]\,d\mu}{\int_{\Omega}u^{4}\,d\mu}\notag\\
&=4\Big(\lambda_{1;f}^2-\frac{2}{3}\lambda_{1;f}^2
+\frac{1}{3}\lambda_{1;f}I_{2,2}\Big)\notag\\
&\leq\frac{16}{3}\lambda_{1;f}^2,
\end{align}
where we used \eqref{0909-Prof-3}.

On the other hand, from the Rayleigh-Ritz inequality, we have
$$\Lambda_{1;f}\leq\frac{\int_{\Omega}(\Delta_f\varphi)^2\,d\mu}{\int_{\Omega}|\nabla\varphi|^2\,d\mu}.$$
Hence,
\begin{align}\label{0916-Prof-4}
\Lambda_{1;f}\lambda_{1;f}&\leq\frac{\int_{\Omega}
(\Delta_f\varphi)^2\,d\mu}{\int_{\Omega}|\nabla\varphi|^2\,d\mu}
\frac{\int_{\Omega}|\nabla\varphi|^2\,d\mu}
{\int_{\Omega}\varphi^2\,d\mu}\notag\\
&=\frac{\int_{\Omega}(\Delta_f\varphi)^2\,d\mu}
{\int_{\Omega}\varphi^2\,d\mu}\notag\notag\\
&\leq\frac{16}{3}\lambda_{1;f}^2,
\end{align}
which shows that $\Lambda_{1;f}\leq \frac{16}{3}\lambda_{1;f}$.

When $\alpha\neq2$, we let $\beta=\frac{1}{2}(\alpha+2)$. Then, if
\begin{align}\label{0916-Prof-5}
{\rm Ric}_f+2(\beta-2)\frac{\nabla^2u}{u}\geq0,
\end{align}
we have
\begin{align}\label{0916-Prof-6}
\frac{\Big[2(2\beta-3)+\frac{4}{3}\Big]
-\frac{4}{3}\sqrt{3(2\beta-3)+1}}
{2(2\beta-3)^2}
&\leq\frac{I_{2,\beta}}{\lambda_{1;f}}\notag\\
&\leq \frac{\Big[2(2\beta-3)+\frac{4}{3}\Big]
+\frac{4}{3}\sqrt{3(2\beta-3)+1}}{2(2\beta-3)^2},
\end{align}
where
$$
I_{2,\beta}=\frac{\int_{\Omega} u^{2\beta-4}|\nabla u|^{4}\,d\mu}{\int_{\Omega} u^{2\beta-2}|\nabla u|^{2} \,d\mu}.
$$
Let $\varphi=u^{\beta}$. Thus, we get
\begin{align}\label{0916-Prof-8}
\Gamma_{1;f}
&\leq \frac{\int_{\Omega}(\Delta_f\varphi)^2\,d\mu}
{\int_{\Omega}\varphi^2\,d\mu}\notag\\
&=\frac{\beta^2\int_{\Omega}\Big[u^{2\beta-2}
(\Delta_fu)^2+2(\beta-1)u^{2\beta-3}|\nabla u|^2\Delta_fu+(\beta-1)^2u^{2\beta-4}|\nabla u|^4\Big]\,d\mu}{\int_{\Omega}u^{2\beta}\,d\mu}\notag\\
&=\beta^2\Big[\lambda_{1;f}^2-\frac{2(\beta-1)}
{2\beta-1}\lambda_{1;f}^2+\frac{(\beta-1)^2}{2\beta-1}
\lambda_{1;f}I_{2,\beta}\Big]\notag\\
&=\frac{\beta^2}{2\beta-1}
[\lambda_{1;f}^2+(\beta-1)^2\lambda_{1;f}I_{2,\beta}].
\end{align}
Combining the above inequalities, we obtain
\begin{align}\label{0916-Prof-9}
\Gamma_{1;f}&\leq \frac{\beta^2}{2\beta-1}\Big[1+\frac{(\beta-1)^2}
{\Big(2\beta-\frac{7}{3}\Big)-\frac{2}{3}\sqrt{2(3\beta-4)}}\Big] \lambda_{1;f}^2\notag\\
&=\frac{\beta^2}{2\beta-1}\Big[1+\frac{1}{3}
\Big(\frac{(3\beta-4)+1}{\sqrt{2(3\beta-4)}-1}\Big)^2\Big]\lambda_{1;f}^2.
\end{align}
It follows from \eqref{0916-Prof-9} by taking $\beta=\frac{4}{3}$ that
\begin{align}\label{0823-Sec-13}
\Gamma_{1;f}\leq&\frac{64}{45}\lambda_{1;f}^2,
\end{align}
and hence $\Lambda_{1;f}\leq \frac{64}{45}\lambda_{1;f}$.
We complete the proof of Theorem \ref{2-thm2}.

\bibliographystyle{Plain}

\end{document}